\documentclass[reqno, 10pt]{amsart}

\usepackage{a4wide, mathrsfs, amssymb}
\usepackage{hyperref}

\numberwithin{equation}{section}

\theoremstyle{plain}
\newtheorem{theorem}{Theorem}[section]

\theoremstyle{definition}

\theoremstyle{remark}

\newcommand\eps{\ensuremath{\varepsilon}}

\newcommand\norm[1]{\ensuremath{\| {#1} \|}}
\newcommand\Wkp{\ensuremath{W^{2, n/2}_{\mathrm{loc}}(M)}}
\newcommand\Ln{\ensuremath{L^{n/2}_{\mathrm{loc}}(M)}}

\newcommand\g{\ensuremath{\mathbf{g}}}
\newcommand\ghat{\ensuremath{\widehat{\mathbf{g}}}}

\newcommand\R{\ensuremath{\mathbb{R}}}

\newcommand\la{\langle}
\newcommand\ra{\rangle}

\newcommand\bel[1]{\begin{equation}\label{#1}}
\newcommand\ee{\end{equation}}

\newcommand\D{\ensuremath{\mathscr{D}}}

\begin{document}
\title{A positive mass theorem for low-regularity metrics}
\date{\today. Preprint UWThPh-2012-18}

\author[J.D.E.\ Grant]{James~D.E.\ Grant}
\address{\href{http://gravity.univie.ac.at}{Gravitationsphysik} \\ Fakult{\"a}t f{\"u}r Physik \\
\href{http://www.univie.ac.at/de/}{Universit{\"a}t Wien} \\ Boltzmanngasse 5 \\ 1090 Wien \\ Austria}
\email{\href{mailto:james.grant@univie.ac.at}{james.grant@univie.ac.at}}
\urladdr{\href{http://jdegrant.wordpress.com}{http://jdegrant.wordpress.com}}

\author[N.\ Tassotti]{Nathalie Tassotti}
\address{\href{http://www.mat.univie.ac.at/home.php}{Fakult\"{a}t f\"{u}r Mathematik} \\
\href{http://www.univie.ac.at}{Universit\"{a}t Wien} \\ Nordbergstra\ss e 15 \\ 1090 Vienna \\ Austria}
\email{\href{nathalie.tassotti@univie.ac.at}{nathalie.tassotti@univie.ac.at}}
\urladdr{\href{http://nathalietassotti.wordpress.com}{http://nathalietassotti.wordpress.com}}


\begin{abstract}
We prove a positive mass theorem for continuous Riemannian metrics in the Sobolev space $W^{2, n/2}_{\mathrm{loc}}(M)$. We argue that this is the largest class of metrics with scalar curvature a positive a.c. measure for which the positive mass theorem may be proved by our methods.
\end{abstract}

\thanks{The work of NT was supported by a Forschungsstipendium from the University of Vienna.}

\maketitle
\thispagestyle{empty}

\section{Statement of results}

Let $(M, \g)$ be a complete, asymptotically flat, smooth Riemannian manifold without boundary of dimension $n \ge 3$. As such, there exists a compact subset $K \subset M$ such that $M \setminus K$ is diffeomorphic to $\cup_{i=1}^N N^i$, with each end $N^i$ diffeomorphic to $\R^n$ minus a compact ball. We assume that the metric $\g$ is smooth on $M \setminus K$, and satisfies the asymptotic conditions required for the validity of the smooth positive mass theorem.%
\footnote{We assume throughout that either $n \le 7$ or that $M$ is a spin-manifold, so that the classical positive mass theorem is valid on $M$.}
We assume that, globally, the metric $\g$ lies in the space $\Wkp \cap C^0(M)$. The Sobolev embedding theorem and continuity of $\g$ implies that the connection lies in $L^n_{\mathrm{loc}}(M)$ and that the scalar curvature, $s_{\g}$, is a well-defined element of $\Ln$. We impose that $s_{\g}$ is non-negative in the distributional sense, i.e.
\[
\la s_{\g}, \varphi^2 \ra \ge 0, \qquad \forall \varphi \in \D(M).
\]
On the set $M \setminus K$, $\g$ is smooth, so this simply states that $s_{\g} \ge 0$. However, inside the set $K$ the metric $\g$ is not assumed $C^2$, and this condition is imposed only in the weak sense.

\medskip
Our main result is the following.

\begin{theorem}
\label{thm}
Let $(M, \g)$ be a Riemannian manifold as above with non-negative scalar curvature in the distributional sense. Then the ADM mass of $(M, \g)$ is non-negative.
\end{theorem}

We outline the proof of this result below. (Full details will appear in~\cite{pmt2}.) It is an adaption of the approach of Miao~\cite{Miao} (see, also,~\cite{McS, Lee}), who considered manifolds that contain hypersurfaces over which the mean curvature could be discontinuous. In~\cite{Miao}, one has $\delta$-function like scalar curvature with singular support on the hypersurface. In our approach, the metric may be non-smooth in an arbitrary compact set, but the non-smooth behaviour is less singular.

\section{Sketch of proof of Theorem~\ref{thm}}

By density of smooth metrics in $\Wkp \cap C^0(M)$, we construct a one-parameter family of smooth metrics $\{ \g_t \}_{t > 0}$ such that $\g_t$ converge to $\g$ both locally uniformly and in $\Wkp$ as $t \to 0$. Enlarging the compact set $K$ if necessary, we may assume that the metrics $\g_t$ coincide with the metric $\g$ on the set $M \setminus K$.%
\footnote{The $\g_t$ may be constructed, for example, as follows. Take a finite collection of charts that cover $K$, and combine these with the ends $N_i$ to give an open cover of $M$. We smooth the metric on the charts covering $K$ by convolution with mollifiers, leaving the metric on the ends unchanged. Adding the resulting metrics together with a smooth partition of unity then yields the smooth approximations that we require.}
Since, by construction, $\g_t = \g$ on $M \setminus K$, the metrics $\g_t$ have the same asymptotic behaviour and mass as $\g$. Moreover, $\g_t$ converges to $\g$ in $L^{\infty}(K)$ and $W^{2, n/2}(K)$ as $t \to 0$. The former fact, along with the property that $\g_t = \g$ on $M \setminus K$ implies that there exists $\rho(t) \ge 1$ with the property that
\bel{equivalence}
\frac{1}{\rho(t)} \g_t \le \g \le \rho(t) \g_t
\ee
as bilinear forms on $M$, with $\rho(t) \to 1$ as $t \to 0$. From our assumption that $\g$ is continuous and lies in $\Wkp$, it follows that $s_{\g_t} \to s_{\g}$ in $\Ln$ as $t \to 0$. In particular, $s_{\g_t} = s_{\g}$ on $M \setminus K$ and $\norm{s_{\g_t} - s_{\g}}_{L^{n/2}(K)} \to 0$ as $t \to 0$. We then have, for all $\varphi \in \D(M)$,
\begin{align*}
\int_M s_{\g_t} \varphi^2 \, d\mu &= \la s_{\g}, \varphi^2 \ra
+ \int_M \left( s_{\g_t} - s_{\g} \right) \varphi^2 \, d\mu
\ge \int_K \left( s_{\g_t} - s_{\g} \right) \varphi^2 \, d\mu
\\
&\ge - \left\Vert s_{\g_t} - s_{\g} \right\Vert_{L^{n/2}(K)} \left\Vert \varphi \right\Vert_{L^{2n/(n-2)}(K)}^2.
\end{align*}
It follows from this inequality that the negative part of the scalar curvature of $\g_t$ satisfies
\bel{smin0}
\left\Vert \left[ s_{\g_t} \right]_- \right\Vert_{L^{n/2}(M)}
\le \left\Vert s_{\g_t} - s_{\g} \right\Vert_{L^{n/2}(M)}
\to 0 \quad \mbox{as $t \to 0$}.
\ee

We recall~\cite[Lemma~3.1]{SY} that, for each metric $\g_t$, there exists a constant $c_1[\g_t] > 0$ such that for any test function $\varphi \in \D(M)$ we have
\[
\norm{\varphi}_{L^{2n/(n-2)}(M, \g_t)}^2 \le c_1[\g_t] \, \norm{\nabla_{\g_t} \varphi}_{L^2(M, \g_t)}^2.
\]
A straightforward calculation using~\eqref{equivalence} shows that the Sobolev constants $c_1[\g_t]$ obey
\bel{Sobequiv}
\frac{1}{\rho(t)^n} \, c_1[\g_t] \le c_1[\g] \le \rho(t)^n \, c_1[\g_t], \qquad t \ge 0.
\ee

Let $c_n = \frac{n-2}{4(n-1)}$. It follows from~\eqref{smin0}, \eqref{Sobequiv} and the fact that $\rho(t) \to 1$ as $t \to 0$ that there exists $t_0 > 0$ such that
\bel{SYcondition}
c_n c_1[\g_t] \, \norm{\left[ s_{\g_t} \right]_-}_{L^{n/2}(M, \g_t)} \le \frac{1}{2}, \qquad t \le t_0.
\ee
An inspection of the proof of Lemma~3.2 in~\cite{SY} (see, also,~\cite[Lemma~4.1]{Miao}) shows that~\eqref{SYcondition} implies that, for each $t \le t_0$, there exists a unique positive solution of the equation
\[
\Delta_{\g_t} u_t + c_n \left[ s_{\g_t} \right]_- u_t = 0
\]
on $M$ with the property that $u_t(x) = \frac{A_t}{r^{n-2}} + \omega_t(r)$ as $r \to \infty$, where $A_t$ is constant, $\omega_t(r) = O(r^{1-n})$ and $\partial\omega_t(r) = O(r^{-n})$ as $r \to \infty$. As in the proof of~\cite[Proposition~4.1]{Miao}, we deduce that the functions $w_t := u_t - 1$ satisfy the inequalities
\begin{align}
\norm{\nabla_{\g_t} w_t}_{L^2(M, \g_t)}^2
\le c_n &\left[ \norm{\left[ s_{\g_t} \right]_-}_{L^{n/2}(M, \g_t)} \norm{w_t}_{L^{2n/(n-2)}(M, \g_t)}^2 \right.
\nonumber
\\
& \hskip 2cm \left. + \norm{\left[ s_{\g_t} \right]_-}_{L^{2n/(n+2)}(M, \g_t)} \norm{w_t}_{L^{2n/(n-2)}(M, \g_t)} \right]
\label{dwbound}
\end{align}
and
\bel{wbound}
\norm{w_t}_{L^{2n/(n-2)}(M, \g_t)} \le 8 c_n^2 c_1[\g_t]^2 \ \norm{\left[ s_{\g_t} \right]_-}_{L^{2n/(n+2)}(M, \g_t)}.
\ee
Equation~\eqref{wbound} implies that $\norm{w_t}_{L^{2n/(n-2)}(M, \g_t)} \to 0$ as $t \to 0$ and therefore, from~\eqref{dwbound}, that $\norm{\nabla_{\g_t} w_t}_{L^2(M, \g_t)}^2 \to 0$ as $t \to 0$.

Finally, we construct conformal rescaled metrics ${\ghat}_t := u_t^{4/(n-2)} \g_t$. These metrics are asymptotically flat and, by construction, have non-negative scalar curvature. Moreover, the mass of $(M, \ghat_t)$ is related to that of $(M, \g_t)$ by the relation~\cite[Eq.~48]{Miao}
\[
m(\g_t) = m(\ghat_t) + \frac{n-1}{(n-2) \omega_{n-1}} \int_M \left[ |\nabla_{\g_t} u_t|_{\g_t}^2 - c_n \left[ s_{\g_t} \right]_- u_t^2 \right] \, d\mu_{\g_t}.
\]
As noted above, the first term in the integral on the right-hand-side converges to $0$ as $t \to 0$. For the second term, we have
\[
\left| \int_M \left[ s_{\g_t} \right]_- u_t^2 \, d\mu_{\g_t} \right|
= \left| \int_K \left[ s_{\g_t} \right]_- u_t^2 \, d\mu_{\g_t} \right|
\le \norm{\left[ s_{\g_t} \right]_-}_{L^{n/2}(K, \g_t)} \norm{u_t}_{L^{2n/(n-2)}(K, \g_t)}^2.
\]
Moreover,
\begin{align*}
\norm{u_t}_{L^{2n/(n-2)}(K, \g_t)}
&\le \norm{1}_{L^{2n/(n-2)}(K, \g_t)} + \norm{w_t}_{L^{2n/(n-2)}(K, \g_t)}
\\
&= |K|_{\g_t}^{\frac{n-2}{2n}} + \norm{w_t}_{L^{2n/(n-2)}(K, \g_t)} \to |K|_{\g}^{\frac{n-2}{2n}} \quad \mbox{as $t \to 0$}.
\end{align*}
From compactness of $K$, we have $|K|_{\g} < \infty$. Since $\norm{\left[ s_{\g_t} \right]_-}_{L^{n/2}(K, \g_t)} \to 0$, it follows that $\int_M \left[ s_{\g_t} \right]_- u_t^2 \, d\mu_{\g_t} \to 0$ as $t \to 0$. Therefore $m(\ghat_t) \to m(\g_t)$ as $t \to 0$. Since $m(\g_t) = m(\g)$, and $m(\ghat_t) \ge 0$ for $t > 0$, we deduce that $m(\g) \ge 0$.

\section{Final remarks}

If the scalar curvature of a metric is well-defined and non-negative as a distribution, it is automatically a positive measure. In~\cite{Miao, McS, Lee}, this measure has a singular part, but the a.c. part is continuous. In our approach the measure is purely a.c., but only lies in $L^{n/2}$ on the compact set $K$. Our class of metrics is therefore \lq\lq complementary\rq\rq\ to that considered in~\cite{Miao, McS, Lee}. For example, the metrics considered in~\cite{Miao, McS, Lee} are Lipschitz, whereas the metrics we consider are not even H\"{o}lder continuous in general. It would be of interest to know whether a combination of our results and those of~\cite{Miao, McS, Lee} can be attained.

We believe that our class of metrics is the largest class of metrics with purely a.c. scalar curvature for which the positive mass theorem holds, at least with this method of proof. If we consider metrics in $\Wkp$ without the continuity assumption, then generally the scalar curvature is not well-defined as a distribution. If we assume that the metric is in $\Wkp$ and locally bounded, then $s_{\g}$ is well-defined as a distribution, but the smoothed metrics $\g_t$ will not converge locally uniformly to $\g$, so we will not have the estimates~\eqref{equivalence} and~\eqref{smin0} that were crucial in our arguments. It therefore seems that continuity of the metric is essential. The fact that we require an $L^{n/2}$ bound on $\left[ s_{\g_t} \right]_-$ in order to apply~\cite[Lemma~3.2]{SY} also suggests that the Sobolev space $\Wkp$ is optimal for our approach. Note that this integral bound on $\left[ s_{\g_t} \right]_-$ arises naturally in our approach rather than the pointwise bound found in~\cite{Miao, McS, Lee}.

On a related topic, there does not appear to be a rigidity result for our class of metrics. It is conceivable that there exist continuous asymptotically flat metrics in $\Wkp$ with non-negative scalar curvature and zero mass that are not flat, although we have been unable to construct an explicit example. The proof of the rigidity property in~\cite{Miao} involves showing that $\sup_M |u_t - 1| \to 0$ as $t \to 0$. One can adapt this argument to metrics that lie in $W^{2, n/2+\eps}_{\mathrm{loc}}(M)$ for some $\eps > 0$. However, the Moser iteration argument required to derive the $L^{\infty}$ bound on $|u_t - 1|$ breaks down when the metric only lies in $\Wkp$. (Generally, one obtains a bound in an appropriate Orlicz space.) It is therefore possible, for example, that the metrics in our class cannot generally be approximated in $C^0$ by smooth asymptotically flat metrics with non-negative scalar curvature.

\medskip

These issues will be reported upon elsewhere~\cite{pmt2}.

\end{document}